\documentclass[11pt, a4paper, oneside]{article}
\usepackage{amsmath, amssymb, amsthm, amsfonts, amsxtra, latexsym, amscd,
pb-diagram,  graphics, hyperref}
\usepackage [all]{xy}

\theoremstyle{plain}

\newcommand{\tx}{\otimes }

\newcommand{\ri}{\rightarrow }

\newcommand{\Fm}{\widetilde{F}}

\DeclareMathOperator{\Coker}{Coker}

\DeclareMathOperator{\Hom}{Hom}

\DeclareMathOperator{\Aut}{Aut}

\DeclareMathOperator{\Ker}{Ker}
\DeclareMathOperator{\Dis}{Dis}

\newtheorem{thm}{\bf Theorem}
\newtheorem{lem}[thm]{\bf Lemma} 
\newtheorem{pro}[thm]{\bf Proposition} 

\theoremstyle{definition}

\begin{document}

\input diagxy

\centerline{\Large  \bf Strict Gr-categories and applications on}
\centerline{\Large  \bf classification of extensions of groups of
the type} \centerline{\Large  \bf of a crossed module}
\vspace{0.5cm}
 \centerline{\bf
Nguyen Tien Quang}
 \centerline{\it Hanoi National University of
Education, Department of Mathematics,}
 \centerline{\it E-mail:
cn.nguyenquang@gmail.com}\vspace{0.25cm} \centerline{\bf Pham Thi
Cuc} \centerline{\it Hongduc University, Natural Science
Department,}
\centerline{\it E-mail: cucphamhd@gmail.com} \vspace{0.25cm}
 \centerline{\bf Nguyen
Thu Thuy} \centerline{\it Academy of Finance, Science Faculty,}
\centerline{\it E-mail: ntthuy11@gmail.com}

\begin{abstract}
In this paper we state some applications of Gr-category theory on
the classification of crossed modules and  on the classification of
extensions of groups of the type of a crossed module.
\end{abstract}

\section{Introduction}
 Theory of  Gr-categories, or 2-groups, and its generality are
getting more and more applications.   The relationship among
Gr-categories, cohomology of groups, and extensions of groups is
stated in \cite{CQT}. The results on group extensions of the type of
a crossed module are also presented by cohomology of groups in \cite
{Br99, Br94, Br96}. These motivate our studies on the representation
of the concepts related to crossed modules by Gr-categories, and
then one can apply the results  of Gr-categories for crossed
modules.

A {\it Gr-category} \cite{Sinh75} is a monoidal category in which
every morphism is invertible and every object has a weak inverse.
(Here, a weak inverse of an object $x$ is an object $y$ such that
$x\otimes y$ and $y\otimes x$ are both isomorphic to the unit
object.)

A {\it strict} Gr-category  is a strict monoidal category in which
every morphism is invertible and every object has a strict inverse
(so that $x\otimes y$ and $y\otimes x$ are actually equal to the
unit object).

Strict Gr-categories can be identified with crossed modules, thus
Gr-categories in general can be seen as a weakening of crossed
modules.

A strict Gr-category is a group object in {\bf Cat}. R. Brown and C.
Spencer   \cite{Br76} studied group objects in the category of
groupoids under the name {\it $\mathcal G$-groupoid}, or {\it
group-groupoid} by R. Brown and O. Mucuk in \cite{Br94}. In
\cite{Br76} the authors stated the construction of a crossed module
from a $\mathcal G$-groupoid and that of a $\mathcal G$-groupoid
from a crossed module. They are repeated by Forrester - Barker
\cite{F-B} as a $\mathcal G$-groupoid is replaced by a strict
Gr-category (see also J. Baez and A. Lauda \cite{Baez}) thanks to
the notion of {\it internal categories}.

In present paper, we show directly the above constructions based on
strict Gr-categories, and prove that the strong homotopy category
$Ho\mathbf{Grstr}$ of strict Gr-categories and single Gr-functors is
isomorphic to the category $\mathbf{Cross}$ of crossed modules
(Theorem \ref{td}).

One obtains the first applications of the Gr-category theory  as
follows. In the classical method, the factor set of a homomorphism
$h:G\ri H$ is only defined if $h$ is a surjection. For the general
case, we can overcome this limitation by building a Gr-functor
between two appropriate Gr-categories. Then the results of
Gr-category theory become the useful device to study crossed
modules.

For any crossed module $B\ri D$ and a group homomorphism $\psi:Q\ri$
Coker$d$ of an extension problem of  type $B\ri D$ (see
\cite{Br94,Br96}), we construct a strict Gr-category $\mathbb P$.
Then, for a Gr-functor $F:\Dis Q\ri \mathbb P$ we can determine an
associated extension of the type of a crossed module. Therefore,
there is a bijection (Theorem \ref{schr})
 $$\mathrm{Hom}_{(\psi,0)}[\Dis Q,\mathbb P]\leftrightarrow \mathrm{Ext}_{B\ri D}(Q,B,\psi),$$
 where $\mathrm{Ext}_{B\ri D}(Q,B,\psi)$ is the set of equivalence classes of  extensions
 of  $B$ by $Q$ of   type $B\ri D$ inducing $\psi$. This result contains The
 classification Theorem of  R. Brown and
  O. Mucuk (Theorem 5.2 \cite{Br94}).

\section{Preliminaries}
For convenience, we recall here some well-known results on
Gr-categories and Gr-functors (see \cite{CQT}).


We often denote by  $(\mathbb G ,\tx, I, \mathbf{a}, \mathbf{l},
\mathbf{r})$  a {\it Gr-category}.   If $(F, \Fm, F_\ast)$ is a
monoidal functor between Gr-categories, it is called a {\it
Gr-functor.} Then the isomorphism $ F_\ast:I'\ri FI$ can be deduced
from $F$ and $\Fm$. Hereafter, we refer to $(F, \Fm)$ as a
Gr-functor.

Two Gr-functors $(F, \Fm)$ and $(F', \Fm')$ from $\mathbb G$ to
$\mathbb G'$ are {\it homotopic} if there is a \emph{natural
monoidal equivalence}, or a {\it homotopy} $\alpha:(F,\Fm,F_\ast)\ri
(F',\Fm',F_\ast')$ which is a natural isomorphism  such that
\begin{equation*}
F'_\ast=\alpha_I\circ F_\ast.
\end{equation*}

A Gr-category is equivalent to one of  type $(\Pi, A),$ which can be
described as follows. The set
$\pi_0\mathbb G$ of isomorphism classes of the objects in $\mathbb
G$ is a group with the operation induced by the tensor product in
$\mathbb G,$ and the set $\pi_1\mathbb G$ of automorphisms of the
unit object $I$ is an abelian group with the operation, denoted by
+, induced by the composition of morphisms. Moreover, $\pi_1\mathbb
G$ is a $\pi_0\mathbb G$-module under the action
\begin{equation*}\label{ct0a}su=\gamma_{X}^{-1}\delta_{X}(u),\ X\in s,\ s\in \pi_0\mathbb G,\ u\in \pi_1\mathbb
G,\end{equation*} where $\delta_X,\ \gamma_X$ are defined by the
following commutative diagrams
$$\begin{CD}
X @>\gamma_X(u) >> X \\
@A \mathbf{l}_X AA     @AA \mathbf{l}_X A \\
I \tx X @>u \tx id >> I \tx X
\end{CD}
\qquad \qquad\qquad
\begin{CD}
X @>\delta_X(u) >> X \\
@A \mathbf{r}_X AA     @AA \mathbf{r}_X A \\
X \tx I @>id \tx u >> X \tx I
\end{CD}$$
\indent The {\it reduced} Gr-category $S_{\mathbb G}$ of  $\mathbb
G$ is a category whose objects are the elements of $\pi_0\mathbb G$
and whose morphisms are automorphisms $(s,u):s\ri s,$ where $s \in
\pi_0\mathbb G,\ u \in \pi_1\mathbb G$. The composition of two
morphisms is induced by the addition in $\pi_1\mathbb G$
$$(s,u).(s,v)=(s,u+v).$$
\indent The category $S_{\mathbb G}$  is equivalent to $\mathbb G$
by canonical equivalences constructed as follows. For each $s=[X]
\in \pi_0\mathbb G,$ choose a representative $X_s\in \mathbb G$ such
that $X_1=I$, and for each $X \in s$ choose an isomorphism $i_X:\
X_s\rightarrow X$ such that $i_{X_s}=id$. The family $(X_s, i_X)$ is
called a {\it stick} of the Gr-category $\mathbb G$ whenever
$$i_{I\otimes X_s}=\mathbf{l}_{X_s},\; i_{X_s\otimes I}=\mathbf{r}_{X_s}.$$
For given stick $(X_s, i_X)$, we obtain two functors
\[\begin{cases}
G:\mathbb G\ri  S_{\mathbb G}\\
G(X)=[ X]=s\\
G(X \stackrel {f}{\ri}Y)=(s,\gamma_{X_s}^{-1}(i_{Y}^{-1}fi_X))
\end{cases}\qquad\qquad
\begin{cases}
H: S_{\mathbb G}\ri  \mathbb G\\
H(s)=X_s\\
H(s,u)=\gamma_{X_s}(u)
\end{cases}\]
\indent Two functors $G$ and $H$ are categorical equivalences by
natural transformations

$$\alpha=(i_{X}): HG\cong id_{\mathbb G};\qquad\ \beta=id:GH\cong id_{ S_{\mathbb G}}.$$

\noindent They are \emph{canonical equivalences}.

With the structure transport  by the quadruple $(G, H, \alpha,
\beta),$ $S_{\mathbb G}$  becomes a Gr-category together with the
following operations
\begin{gather*}
s\tx t=s.t,\quad s,t\in\pi_0\mathbb G,\\
(s,u)\tx(t,v)=(st,u+sv),\quad u,v\in\pi_1\mathbb G.
\end{gather*}
 \indent The unit constraints of the
Gr-category $S_\mathbb G$ are therefore strict, and its
associativity constraint  is ${\bf a}_{s,r,t}=(srt,k(s,r,t)), $
where
   $$k\in Z^{3}(\pi_0\mathbb G,\pi_1\mathbb G).$$
Moreover, the equivalences $G$ and $ H$ become Gr-equivalences
together with natural isomorphisms
\begin{equation}\label{ct1a}\widetilde{G}_{X,Y}=G(i_X\tx i_Y)\ ,\
 \widetilde{H}_{s,t}=i_{X_s\tx X_t}^{-1}:X_sX_t\ri
 X_{st}.\end{equation}


The Gr-category $S_{\mathbb G}$  is called a \emph{reduction} of the
Gr-category $\mathbb G.$ $S_{\mathbb G}$  is said to be of {\it
type} $(\Pi, A, k)$, or simply  {\it type} $(\Pi,A)$ if
$\pi_0\mathbb G, \pi_1\mathbb G$ are replaced with the group $\Pi$
and the $\Pi$-module $A$, respectively.

Let $\mathbb S=(\Pi, A,k),\ \mathbb S'=(\Pi',A',k')$ be
Gr-categories. A functor $F: \mathbb S\rightarrow \mathbb S'$ is  of
{\it type} $(\varphi,f)$ if
\begin{equation*}
F(x)=\varphi(x),\ \ F(x,u)=(\varphi(x), f(u)),
\end{equation*}
where $\varphi:\Pi\ri \Pi'$, $f:A\ri A'$ are group homomorphisms
satisfying $f(xa)=\varphi (x)f(a)$ for $x\in \Pi, a\in A.$
\begin{lem}[\cite{CQT}\label{bd0}]
Each Gr-functor $(F, \Fm): \mathbb G\rightarrow \mathbb G'$ induces
a Gr-functor $S_F:S_{\mathbb G}\rightarrow S_{\mathbb G'}$ of  type
$(\varphi, f)$, where $\varphi=F_0, f=F_1$ which are determined by
\begin{gather*}
 F_0: \pi_{0}\mathbb G\rightarrow \pi_{0}\mathbb G^{'}, \ \ [ X]\mapsto [FX],\\
F_1:\pi_{1}\mathbb G\rightarrow \pi_{1}\mathbb G^{'}, \ \ u\mapsto
\gamma^{-1}_{FI}(Fu)
\end{gather*}
 Moreover,
$$S_F=G' FH,$$
where $H, G'$ are canonical equivalences.
\end{lem}
Note that if $\Pi'$-module $A'$ is considered as a $\Pi$-module
under the action  $xa'=\varphi(x).a'$, then $f:A\ri A'$ is a
homomorphism of $\Pi$-modules.
 The compatibility of $(F,
\Fm)$ with the associativity constraint gives
$$\varphi^\ast k'-f_\ast k=\partial(g_F),$$
where
\begin{gather*}
 (f_{\ast}k)(x,y,z)=f(k(x,y,z)),\\
(\varphi^{\ast}k')(x,y,z)=k'(\varphi x, \varphi y,\varphi z),
\end{gather*}
and  $g_F:\Pi^2\ri A'$  is a function {\it associated}  to
$\widetilde{F}$.

 If $F:\mathbb S\rightarrow \mathbb S'$
is a functor of  type $(\varphi, f)$, then the function
 \begin{equation}\label{ct1b}
 \xi=\varphi^{\ast}k'-f_{\ast}k
 \end{equation}
is called an {\it obstruction} of the functor $F.$ We have
\begin{pro}[\cite{CQT}\label{md0}]
The functor  $F:\mathbb S\ri \mathbb S' $  of  type  $(\varphi, f)$
induces  a
   Gr-functor if
and only if its obstruction
  $\overline{\xi}$ vanishes
 in  $H^3(\Pi, A')$. Then, there is a bijection
 $$\mathrm{Hom}_{(\varphi,f)}[\mathbb S,\mathbb S']\leftrightarrow
 H^2(\Pi,A'),$$
 where $\mathrm{Hom}_{(\varphi,f)}[\mathbb S,\mathbb S']$ is the set
 of homotopy classes of the Gr-functors of  type
  $(\varphi, f)$ from $\mathbb
 S$ to $\mathbb S'$.
 \end{pro}



\section{Gr-categories associated to a crossed module}
{\bf Definition.} A {\it crossed module} is a quadruple
$(B,D,d,\theta)$ where $d:B\ri D,\;\theta:D\ri$ Aut$B$ are group
homomorphisms making the following diagram commute
\begin{equation*}\begin{diagram}\label{ct4a}
\node{B}\arrow{se,b}{\mu}\arrow[2]{e,t}{d}\node[2]{D}\arrow{sw,b}{\theta}\\
\node[2]{\mathrm{Aut}B}
\end{diagram}\end{equation*}
and satisfying the relation
 \begin{equation}\label{ct4b}d(\theta_x(b))=\mu_x(d(b)),\;x\in D,
b\in B,\end{equation} where $\mu_x$ is an inner automorphism given
by conjugation of $x$.

In present paper, the crossed module $(B,D,d,\theta)$ is sometimes
denoted by
 $B\stackrel{d}{\rightarrow}D$, or $B\rightarrow D$.


For convenience, we denote by the addition for the operation in $B$
and by the multiplication for that in $D$.

The following properties follow from the definition of a crossed
module.

\begin{pro} \label{md2}
Let $(B,D,d,\theta)$ be a crossed module. \\
$\mathrm{i)}\; \mathrm{Ker}d\subset Z(B)$.\\
$\mathrm{ii)} \;\mathrm{Im}d$ is a normal subgroup in $D$.\\
$\mathrm{iii)}$ The homomorphism $\theta$ induces a homomorphism
$\varphi:D\ri \mathrm{Aut}(\mathrm{Ker}d)$ given by
$$\varphi_x=\theta_x|_{\mathrm{Ker}d}.$$
$\mathrm{iv)}\; \mathrm{Ker}d$ is a left $\mathrm{Coker}d$-module
under the action
$$sa=\varphi_x(a),\ \ a\in\mathrm{ Ker}d, \ x\in s\in \mathrm{Coker}d.$$
\end{pro}
 For any crossed module $(B,D,d, \theta)$ we can
construct a strict Gr-category $\mathbb P_{B\rightarrow D}=\mathbb
P$ called the Gr-category {\it associated to} the crossed module
$B\ri D$, as follows.
$$\mathrm{Ob}(\mathbb P)=D, \mathrm{Hom}(x,y)=\{b\in B/x=d(b)y\},$$ where $x,y$ are objects of $\mathbb P$.
 The composition of two
morphisms is given by
$$(x\stackrel{b}{\ri}y\stackrel{c}{\ri}z)=(x\stackrel{b+c}{\ri}z).$$
The tensor operation on objects is given by the multiplication in
the group $D$, and for two morphisms $(x\stackrel{b}{\ri}y),
(x'\stackrel{b'}{\ri}y')$ then
\begin{equation}\label{tx0}
(x\stackrel{b}{\ri}y)\otimes(x'\stackrel{b'}{\ri}y')=(xx'\stackrel{b+\theta_yb'}{\longrightarrow}yy').
\end{equation}
By the definition of a crossed module, we can easily check that
 $\mathbb P$ is a  Gr-category with
the identity constraints.

Conversely, for a strict Gr-category $(\mathbb P, \otimes)$ we
define an {\it associated} crossed module $C_\mathbb
P=(B,D,d,\theta)$ as follows. Set
\begin{equation} D= \mathrm{Ob}(\mathbb P),\notag\end{equation}
\begin{equation}B=\{x\xrightarrow{b}1 / x\in D  \}. \notag \end{equation}
The operations of $D$ and of $B$ are given by
$$xy=x\otimes y,\;\;b+c=b\otimes c, $$
respectively. Then $D$ becomes a group in which the unit is $1$, the
inverse of $x$ is $x^{-1}$ ($x\otimes x^{-1}=1$). $B$ is a group in
which the unit is the morphism   $(1\xrightarrow{id_1} 1)$ and the
inverse of $(x\xrightarrow{b} 1)$ is the morphism
$(x^{-1}\xrightarrow{\overline{b}} 1) \;(b\otimes
\overline{b}=id_1)$.

The homomorphisms $d:B\rightarrow D$ and $\theta: D\rightarrow
\mathrm{Aut} B$ are given by
$$d(x\xrightarrow{b}  1)=x,$$
$$\theta _y(x\xrightarrow{b}1)= (yxy^{-1}\xrightarrow{id_y + b + id_{y^{-1}}}  1),$$ 
respectively. It is easy to see that  $(B,D,d,\theta )$ is a crossed
module.

\noindent{\bf Definition.} A {\it morphism}
$(f_1,f_0):(B,D,d,\theta)\ri (B',D',d',\theta')$ of crossed modules
consists of homomorphisms of groups $f_1:B\ri B'$, $f_0:D\ri D'$
such that the following diagram commutes
\begin{equation}\label{gr1}\begin{diagram}
\node{B}\arrow{e,t}{d}\arrow{s,l}{f_1}\node{D}\arrow{s,r}{f_0}\\
\node{B'}\arrow{e,b}{d'}\node{D'}
\end{diagram}\end{equation}
and $f_1$ is an {\it operator homomorphism}, i.e. for all $x\in D,b\in B,$
\begin{equation}\label{gr2}
f_1(\theta_xb)=\theta'_{f_0(x)}f_1(b).
\end{equation}


\begin{pro}\label{t1}
Let $(f_1,f_0): (B,D,d,\theta)\ri(B',D',d',\theta')$ be a morphism
of crossed modules. Let $\mathbb P$ and $\mathbb P'$ be two
Gr-categories  associated to  crossed modules  $(B,D,d,\theta)$ and
$(B',D',d',\theta')$, respectively.

\noindent $\mathrm{i)}$ There exists a functor $F:\mathbb
P\ri\mathbb P'$  defined by
$F(x)=f_0(x),\;F(b)=f_1(b),$ for $x\in \mathrm{Ob}(\mathbb P)$,
$b\in \mathrm{Mor}(\mathbb P)$.

\noindent $\mathrm{ii)}$ A natural isomorphism
$\widetilde{F}_{x,y}:F(x)F(y)\ri F(xy)$ together with $F$ is a
Gr-functor if and only if $\widetilde{F}_{x,y}$ is a constant
$\widetilde{c}\in\mathrm{Ker}d'$ and
\begin{equation}\theta'_{f_0(x)}(\widetilde{c})=\widetilde{c}.\label{htbs}\end{equation}
\end{pro}
We say that $F$ is a Gr-functor of  form $(f_1,f_0)$.

\begin{proof}
i) By the determination of the Gr-category associated to a crossed
module,  one can easily check that $F$ is a functor.

\noindent ii) We define a natural isomorphism
$$\widetilde{F}_{x,y}:F(x)F(y)\ri F(xy)$$
such that $F=(F,\widetilde{F})$ becomes a Gr-functor. We first see
that if $(x\stackrel{b}{\ri}y)$ and $(x'\stackrel{b'}{\ri}y')$ are
morphisms in $\mathbb P,$ then
\[F(b\otimes
b')=F(xx'\stackrel{b+\theta_yb'}{\longrightarrow}yy')=\big(f_0(xx')\stackrel{f_1(b+\theta_yb')}{\longrightarrow}f_0(yy')\big),\]
\[\begin{aligned}F(b)\otimes F(b')&=\big(f_0(x)\stackrel{f_1(b)}{\longrightarrow}f_0(y)\big)\otimes
\big(f_0(x')\stackrel{f_1(b')}{\longrightarrow}f_0(y')\big)\\
\;&=\big(f_0(x)f_0(x')\xrightarrow{f_1(b)+\theta'_{f_0(y)}f_1(b')}f_0(y)f_0(y')\big).\end{aligned}\]
 Since $f_1$ is a group homomorphism and by
 (\ref{gr2}) we have
\begin{equation}\label{tx1}F(b\otimes b')=F(b)\otimes F(b').\end{equation}

Besides, $F(x)F(y)=F(xy)$ and
$F(x)F(y)=d'(\widetilde{F}_{x,y})F(xy)$ lead to
$d'(\widetilde{F}_{x,y})=1.$ Therefore,
$$\widetilde{F}_{x,y}\in \mathrm{Ker}d'\subset Z(B').$$
The natural isomorphisms $\widetilde{F}_{x,y}$ (if they exist) must
satisfy the diagram \begin{equation}\begin{diagram}
\node{F(x)F(x')}\arrow{e,t}{\widetilde{F}_{x,x'}}\arrow{s,l}{F(b)\otimes
F(b')}\node{F(xx')}\arrow{s,r}{F(b\otimes b')}\\
\node{F(y)F(y')}\arrow{e,b}{\widetilde{F}_{y,y'}}\node{F(yy')}
\end{diagram}\label{bdtn}\end{equation}
%
From (\ref{tx1}) and $\widetilde{F}_{x,y}\in Z(B'),$ this
commutative diagram implies
%
%
 $\widetilde{F}_{x,x'}=\widetilde{F}_{y,y'}.$
%
%
%
Write $\widetilde{F}_{x,y}=\widetilde{c}$. The compatibility of
$(F,\widetilde{F})$ with the associativity constraints implies
$\theta'_{F(x)}(\widetilde{c})=\widetilde{c}$, for $x\in D$.
\end{proof}

A Gr-functor $(F,\widetilde{F}):\mathbb P\ri\mathbb P'$ is {\it
single} if $\widetilde{F}_{x,y}=\widetilde{c}\in \mathrm{Aut}(1')$
for all $x,y\in$ Ob$(\mathbb P)$.

 We state the converse proposition of
Proposition \ref{t1}.
\begin{pro}\label{n1}
Let $\mathbb P$, $\mathbb P'$ be corresponding  Gr-categories
associated to the crossed modules $(B,D,d,\theta)$,
$(B',D',d',\theta')$, and $F:\mathbb P\ri \mathbb P'$ be a single
Gr-functor. Then, $F$ induces a morphism of crossed modules
$(f_1,f_0):(B\ri D)\rightarrow (B'\ri D')$, where $$f_1(b)=F(b),\;\;
f_0(x)=F(x),$$ for all $b\in B, x\in D$.
\end{pro}

\begin{proof}
 Any
$b\in B$  can be considered as a morphism $(db\stackrel{b}{\ri} 1)$
in $\mathbb P$, hence $(F(db)\stackrel{F(b)}{\ri} 1')$ is a morphism
in $\mathbb P'.$

 For each $x,y\in D$, since
$(FxFy\stackrel{\widetilde{F}_{x,y}}{\longrightarrow}F(xy))$ is a
morphism in $\mathbb P',$ and since $\widetilde{F}_{x,y}$ is in
Aut$(1')=$ Ker$d'$, we obtain
$$FxFy=d'(\widetilde{F}_{x,y})F(xy)=F(xy)$$
or $f_0(xy)=f_0(x)f_0(y)$.

For any morphism $(x\stackrel{b}{\rightarrow}y)$ in $\mathbb P,$ we
have $x=d(b)y.$ It follows that
\[
f_0(x)=f_0(d(b)y)=f_0(d(b))f_0(y).\] Besides,
$(f_0(x)\stackrel{f_1(b)}{\rightarrow}f_0(y))$ is a morphism in
$\mathbb P'$, so
\[f_0(x)=d'(f_1(b))f_0(y).\]
Therefore, $f_0(d(b))=d'(f_1(b)),$
for all $b\in B.$ That means the diagram (\ref{gr1}) commutes.

Now, since $\widetilde{F}_{x,y}\in \mathrm{Ker}d'$, the commutative
(\ref{bdtn}) implies (\ref{tx1}).
According to the determination of $F(b)\otimes F(b')$ and
$F(b\otimes b')$ in the proof of Proposition \ref{t1}, we have
\[
f_1(b)+f_1(\theta_yb')=f_1(b)+\theta'_{f_0(y)}f_1(b').\] Hence
$f_1(\theta_yb')=\theta'_{f_0(y)}f_1(b'),$ for all $b'\in B,y\in D,$
or (\ref{gr2}) holds.
\end{proof}

We write
 ${\bf Grstr}$
for the  category of strict Gr-categories and simple Gr-functors.
Two Gr-functors $(F,\widetilde{F})$ and $(G,\widetilde{G})$ are {\it
strong homotopic} if they are homotopic and $F=G$.
 We can define the {\it strong homotopy category}
$Ho\mathbf{Grstr}$ of ${\bf Grstr}$ to be the quotient category with
the same objects, but morphisms are homotopy classes of monoidal
functors. We write
 ${\mathrm{Hom}}_{{\bf Grstr}}[\mathbb P,\mathbb P'] $ for the homsets of the homotopy category, that is,
$${\mathrm{Hom}}_{{\bf Grstr}}[\mathbb P,\mathbb P']=\frac{{\mathrm{Hom}}_{{\bf Grstr}}(\mathbb P,\mathbb P')}
{\text{strong homotopies}}$$
Let {\bf Cross} for the category of crossed modules, we obtain the
similar result to  Theorem 1 \cite{Br76}
\begin{thm}[Classification Theorem]\label{td} There exists an
isomorphism
\[\begin{matrix}
 \Phi:&{\bf Cross} &\ri&Ho{\bf Grstr}\\
&(B\ri D)&\mapsto&\mathbb{P}_{B\ri D}\\
&(f_1,f_0)&\mapsto&[F]
\end{matrix}\]
where $F(x)=f_0(x),F(b)=f_1(b)$, for $x\in \mathrm{Ob}(\mathbb P),
b\in \mathrm{Mor} (\mathbb P)$.
\end{thm}
\begin{proof} Let $\mathbb P,\mathbb P'$ be corresponding Gr-categories associated to
the crossed modules $B\ri D,B'\ri D'$, and $(F,\widetilde{F}),
(G,\widetilde{G}):\mathbb P\ri\mathbb P'$ be two Gr-functors of the
form $(f_1,f_0)$. Then $F=G$ and by Proposition \ref{t1},
$\widetilde{F}=c$ and $\widetilde{G}=c'$ are constants satisfying
the rule (\ref{htbs})
$$\theta'_{f_0(x)}(c)=c, \theta'_{f_0(x)}(c')=c'.$$
 We
show that $a=c-c'$ is a homotopy of $(F,\widetilde{F})$ and
$(G,\widetilde{G})$, i.e., $[F]=[G]$ and hence, $\Phi$ is a map on
the homsets
$$\Phi:\mathrm{Hom}_{{\bf Cross}}(B\ri D, B'\ri D')\rightarrow \mathrm{Hom}_{{\bf Grstr}}
[\mathbb P_{B\ri D},\mathbb P'_{B'\ri D'}].$$  Indeed, the
naturality of $a$ is trivial. The commutativity of the diagram
\[\begin{diagram}
\node{F(x)F(y)}\arrow{e,t}{c}\arrow{s,l}{a\otimes
a}\node{F(xy)}\arrow{s,r}{a}\\
\node{G(x)G(y)}\arrow{e,b}{c'}\node{G(xy)}
\end{diagram}
\]
means that
$$c+a=a+\theta'_{G(x)}(a)++c'.$$
Since $c,c'$ satisfy the rule (\ref{htbs}), so does $a$. Now, since
$c,c'\in\mathrm{Ker}d'$, so the above relation holds.

Since the homotopy between Gr-functors is strong, $\Phi$ is an
injection. By Proposition
 \ref{n1}, every single Gr-functor $F:\mathbb P_{B\ri D}\ri\mathbb P'_{B'\ri
 D'}$ determines a morphism $(f_1,f_0)$ of crossed modules, and
 clearly $\Phi(f_1,f_0)=[F]$. Thus, $\Phi$ is a surjection on morphisms.

By the construction of the Gr-category associated to a crossed
module, $\Phi$ is a bijection on objects. Thus, $\Phi$ is an
isomorphism.
\end{proof}

\section {Extensions of groups of the type of a crossed module}
We now recall the notion of an {\it extension of  groups of the type
of a crossed module} due to Taylor \cite{Tay} and Dedeker \cite{Ded}
(see also \cite{Br96}).

Note that if $B$ is a normal subgroup in $D$, then the system $(B,D,
d, \theta)$ is a crossed module in which $d:B\ri D$ is an inclusion,
$\theta:D\ri \mathrm{Aut}B$ given by conjugation.

\noindent\textbf{Definition.} Let $d: B\rightarrow D$ be a crossed
module. A group \emph{extension} of $B$ by $Q$, \emph{of  type } $d:
B\rightarrow D$ is  a diagram of homomorphisms of groups
\begin{align*}  \begin{diagram}
\xymatrix{ \mathcal E:\;\;\;0 \ar[r]& B \ar[r]^j \ar@{=}[d] &E  \ar[r]^p \ar[d]^\varepsilon & Q \ar[r]& 1, \\
& B \ar[r]^d & D}
\end{diagram}
\end{align*}
where the top row is exact, and the family $(B,E,j,\theta')$ is a
crossed module where $\theta'$ is given by conjugation, and
$(id,\varepsilon)$ is a morphism of crossed modules.

\vspace{10pt} Two extensions of $B$ by $Q$ of  type
$B\xrightarrow{d}D$ are said to be \emph{equivalent } if there is a
morphism of exact sequences
\begin{align*}  \begin{diagram}
\xymatrix{ 0 \ar[r]& B \ar[r]^j \ar@{=}[d] &E \ar[r]^p \ar[d]^\eta &
Q \ar[r]\ar@{=}[d] & 1,&\;\;\;\;E
\ar[r]^\varepsilon&D \\
 0 \ar[r]& B \ar[r]^{j'}  &E'  \ar[r]^{p'}   & Q \ar[r]&
 1,&\;\;\;\;E'
\ar[r]^{\varepsilon'}&D}
\end{diagram}
\end{align*}
such that $\varepsilon'\eta=\varepsilon$. Obviously, $\eta $ is
an isomorphism.


In the diagram
\begin{align} \label{bd10} \begin{diagram}
\xymatrix{ 0 \ar[r]& B \ar[r]^j \ar@{=}[d] &E  \ar[r]^p \ar[d]^\varepsilon & Q \ar[r]\ar@{.>}[d]^\psi & 1, \\
 & B \ar[r]^d  &D  \ar[r]^q   & \text{Coker}d }
\end{diagram}
\end{align}
where $q$ is a canonical homomorphism, since the top row is exact
and $ q\circ \varepsilon \circ j = q\circ d =0,$ there  is a
homomorphism $\psi: Q\rightarrow \text{Coker}d $ such that the right
hand side square commutes. Moreover, $\psi$ depends only on the
equivalence class of the extension.
Our objective is to study the set
$$\mathrm{Ext}_{B\ri D}(Q,B,\psi)$$
 of equivalence
classes of extensions of $B$ by $Q$ of type $B\ri D$ inducing $\psi:
Q\ri \mathrm{Coker} d$. Such group extensions  have been classified
by R. Brown and O. Mucuk (Theorem 5.2 \cite{Br94}).

In the present paper, we use the obstruction theory of Gr-functors
to prove Theorem 5.2 \cite{Br94}. Further, the second assertion of
this theorem can be  seen as a consequence of Schreier Theory
(Theorem \ref{schr}) by Gr-functors between strict Gr-categories
$\mathbb P_{B\ri D}$ and Dis$Q$, where Dis$Q$ is a Gr-category of
type $(Q,0,0)$.

Let $\mathbb P=\mathbb P_{B\ri D}$ be the  Gr-category associated to
crossed module $B\stackrel{d}{\ri}D$. Since $\pi_0\mathbb
P=\mathrm{Coker}d, \;\pi_1\mathbb P=\mathrm{Ker}d$, then
$$S_\mathbb P=(\mathrm{Coker}d, \mathrm{Ker}d,k),\overline{k}\in H^3(\mathrm{Coker}d, \mathrm{Ker}d).$$
The homomorphism $\psi:Q\rightarrow \mathrm{Coker} d$ induces an
{\it obstruction}
$$\psi^*k\in Z^3(Q,\mathrm{Ker}d).$$
Under this notion of obstruction, we state and prove the following
theorem.


\begin{thm}\label{BM}  Let $(B,D,d,\theta)$ be a crossed module and
 $\psi: Q\rightarrow \mathrm{Coker} d$ be a group homomorphism.
 Then,
the vanishing of $\overline{\psi^*k}$ in $ H^3_\psi(Q,\mathrm{Ker}d)
$ is  necessary and sufficient for there to exist an extension of
$B$ by $Q$ of  type $B\rightarrow D$ inducing $\psi$. Further, if
$\overline{\psi^*k} $ vanishes, then the equivalence classes of such
extensions are bijective with $H^2_\psi(Q,\mathrm{Ker}d)$.
\end{thm}


The first assertion of Theorem \ref{BM} follows from the following
lemma.

\begin{lem}\label{mrlk}
For each Gr-functor $(F,\widetilde{F}):\mathrm{Dis} Q\ri \mathbb P$,
 there
exists an extension $\mathcal E_F$ of $B$ by $Q$ of type $B\ri D$
inducing $\psi: Q\ri \mathrm{Coker} d$.
\end{lem}
Such extension $\mathcal E_F$ is called an {\it extension associated
} to the Gr-functor $F$.
\begin{proof} By Lemma \ref{bd0}, $(F,\widetilde{F})$ induces a Gr-functor $S_F:\mathrm{Dis}Q\ri S_\mathbb P$
of type $(\psi,0)$.
 Let $(H,\widetilde{H}):S_\mathbb P\rightarrow \mathbb P$ be a canonical
Gr-functor defined by the stick $(x_s, i_x)$. By (\ref{ct1a}),
\begin{equation*}\label{hn} H(s)= x_s,\; H(s,b)=b,\;\widetilde{H}_{s,r}=-i_{x_s\cdot x_r}.
\end{equation*}

Also by Lemma \ref{bd0}, $(F,\widetilde{F})$ is homotopic to the
composition
$$ \text{Dis} Q\xrightarrow{S_F}  S_\mathbb P\xrightarrow{H}\mathbb P.$$
 So we can choose
$(F,\widetilde{F})=(H,\widetilde{H})\circ(S_F,\widetilde{S_F})$.
Then, by the determination of  $\widetilde{HS_F}$,

\begin{equation}\label{ht12} F(u) = x_s,
\end{equation}
\begin{equation}\label{fn}\widetilde{F}_{u,v} =f(u,v)= h(u,v)-i_{x_s\cdot  x_r}\in B, \end{equation}
for $ u,v\in Q, s=\psi(u),\ r= \psi(v),
h(u,v)=(\widetilde{S_F})_{u,v}.$

By the compatibility of $(F,\widetilde{F})$ with the strict
constraints of Dis$Q$ and $\mathbb P$, $f$ is a ``normal" function
satisfying
\begin{equation} \theta _{Fu}[f(v,t)]+ f(u,vt)=f(u,v)+ f(uv,t). \label{eq1}
\end{equation}

The function $\varphi : Q\rightarrow \Aut B$ given by
\begin{equation*}\label{vp}\varphi (u)= \theta _{Fu}
= \theta _{x_s}\;\; (s=\psi(u))
\end{equation*}
satisfies the rule \begin{equation}\label{eq2}\varphi (u)\varphi
(v)=\mu [f(u,v)] \varphi (uv).\end{equation}

 Indeed, by (\ref{fn})
and $h(u,v)\in \Ker d$, we have
\begin{align*} \varphi (u)\varphi (v)&= \theta_{x_s}\cdot \theta_{x_r} = \theta _{x_s\cdot x_r} \notag{}\\
&= \theta [d(-i_{x_s\cdot x_r})x_{sr}]
= \theta [d(-i_{x_s\cdot x_r})]\theta (x_{sr})\notag{}\\
&= \mu(-i_{x_s\cdot x_r}) \varphi (uv)=\mu [f(u,v)] \varphi (uv).
\end{align*}
The pair $(\varphi , f)$ satisfies \eqref{eq1} and \eqref{eq2}, so
it is a {\it factor set}, and therefore one can define the crossed
product $E_0=[B,\varphi ,f, Q]$ (see \cite{MacL63}), i.e.,
$ E_0=B\times Q$ with the operation
\begin{equation*}
(b,u)+(b',u')= (b+\varphi (u)b' + f(u,u'),uu').
 \end{equation*}
 The unit is $(0,1)$ and $-(b,u)=(c,u^{-1})$, where
$\varphi(u)c=-b-f(u,u^{-1})$. Then, we have an exact sequence
$$\mathcal E_F: 0\ri B\stackrel{j_0}{\ri}E_0\stackrel{p_0}{\ri}Q\ri
1,$$ where
$$ j_0(b)=(b,1);\  p_0(b,u)= u, \ b\in B, u\in Q.$$
Clearly, $j_0(B)$ is a normal subgroup in $E_0$, so $j_0:B\ri E_0$
is a crossed module in which the action $\theta':E_0\ri \mathrm{Aut}
B$ is given by conjugation.


We define a group homomorphism $\varepsilon: E_0\ri D$ by
$$\varepsilon(b,u)=d(b)x_{\psi(u)},\;(b,u)\in E_0,$$
where $x_{\psi(u)}$ is a representative of $\psi(u)$ in $D$. Then,
$(id,\varepsilon)$ is a morphism of crossed modules. Indeed,
$\varepsilon\circ j_0=d$. Further, for all $(b,u)\in E_0, c\in B$
$$\theta'(b,u)(c)=j_0^{-1}[\mu_{(b,u)}(c,0)]=\mu_bc+\varphi(u)c,$$
$$\theta_{\varepsilon(b,u)}c=\theta_{dbx_{\psi(u)}}c=\mu_bc+\varphi(u)c.$$
Hence, $\theta'(b,u)(c)=\theta_{\varepsilon(b,u)}c$.


Therefore, we have an extension
\begin{equation*}  \begin{diagram}
\xymatrix{\mathcal E_F:\;\;\;\;\; 0 \ar[r]& B \ar[r]^{j_0} \ar@{=}[d] &E_0  \ar[r]^{p_0} \ar[d]^\varepsilon & Q \ar[r]\ar@{.>}[d]^\psi & 1, \\
 & B \ar[r]^d  &D  \ar[r]^q   & \text{Coker}d }
\end{diagram}
\end{equation*}

For all $u\in Q$,  $q\varepsilon(0,u)= q(x_{\psi(u)})= \psi(u),$
i.e., $\mathcal E_F$ induces $\psi: Q\rightarrow \text{Coker}d$.
\end{proof}

{\bf The proof of Theorem \ref{BM}}
\begin{proof}
Let's recall that $\mathbb P$ is the  Gr-category associated to the
crossed module $B\ri D$. Then, its reduced Gr-category is $S_\mathbb
P=(\mathrm{Coker}d,\mathrm{Ker}d,k)$, where $k\in
Z^3(\mathrm{Coker}d,\mathrm{Ker}d)$.  The pair
$$(\psi,0):(Q,0,0)\ri (\Coker d, \Ker d,0)$$
has $-\psi^*k$ as an obstruction. By the assumption,
$\overline{\psi^*k}=0$, hence by Proposition \ref{md0} it determines
a Gr-functor $(\Psi, \widetilde{\Psi}):\mathrm{Dis} Q\ri S_\mathbb
P$. Then the composition of $(\Psi, \widetilde{\Psi})$ and
$(H,\widetilde{H}): S_\mathbb P{\ri}\mathbb P$ is a Gr-functor
$(F,\widetilde{F}):\mathrm{Dis} Q{\ri} \mathbb P$, and by Lemma
\ref{mrlk} we obtain an associated extension $\mathcal E_F$.

Conversely, suppose that there is an extension as in the diagram
(\ref{bd10}). Let $\mathbb P'$ be the category associated to the
crossed module $B\ri E$. By Proposition \ref{t1}, there is a
Gr-functor $F:\mathbb P'\ri\mathbb P.$ Since the reduced Gr-category
of $\mathbb P'$ is Dis$Q$, by Lemma \ref{bd0} $F$ induces a
Gr-functor of type $(\psi,0)$ from Dis$Q$ to
$(\mathrm{Coker}d,\mathrm{Ker}d,k)$. Now, by Proposition \ref{md0},
the obstruction of the pair $(\psi,0)$  must vanish in
$H^3(Q,\mathrm{Ker}d)$, i.e. $\overline{\psi^*k}=0$.

The final assertion  of the theorem is obtained  from Proposition
\ref{t1}.
\end{proof}


\begin{thm}[Schreier  Theory for extensions of the type of a crossed module]\label{schr}  There is a bijection
$$
\Omega: \mathrm{Hom}_{(\psi,0)}[\mathrm{Dis}Q,\mathbb
P]\rightarrow\mathrm{Ext}_{B\ri D}(Q, B,\psi).$$
\end{thm}
\begin{proof}

{\it Step 1: Gr-functors $(F,\widetilde{F}) $ and
$(F',\widetilde{F'}) $ are homotopic if and only if the
corresponding associated extensions $\mathcal E_F, \mathcal E_{F'}$
are equivalent.}




Let $F, F':$ Dis$Q\ri\mathbb P$ be homotopic
by a homotopy $\alpha:F\ri F'$. By Lemma \ref{mrlk}, there exist the
extensions $\mathcal E_F$ and $\mathcal E_{F'}$ associated to $F$
and $ F'$, respectively. By the definition of Gr-morphisms, the
following diagram commutes
\[\begin{diagram}
\node{FuFv}\arrow{e,t}{\widetilde{F}_{u,v}}\arrow{s,l}{\alpha_u\otimes
\alpha_v}\node{F(uv)}\arrow{s,r}{\alpha_{uv}}\\
\node{F'uF'v}\arrow{e,b}{\widetilde{F'}_{u,v}}\node{F'(uv)}
\end{diagram}\]
That means $$\widetilde{F}_{u,v}+\alpha_{uv}=\alpha_u\otimes
\alpha_v+\widetilde{F}'_{u,v}$$ \noindent (the sum in $B$). By
relation (\ref{tx0}),
\begin{equation}\label{eq5}
f(u,v)+\alpha_{uv}=\alpha_u+\theta_{F'u}(\alpha_v)+f'(u,v),
\end{equation}
where $f(u,v)=\widetilde{F}_{u,v}, f'(u,v)=\widetilde{F}'_{u,v}$.
Now we
 set
\begin{align*} \alpha^*: E_F&\rightarrow E_{F'},\\
(b,u)&\mapsto (b+\alpha_u,u). \end{align*}

Then $\alpha^*$ is a homomorphism thanks to the relation
(\ref{eq5}). Further, the following diagram commutes
\begin{align*}  \begin{diagram}
\xymatrix{ 0 \ar[r]& B \ar[r]^j \ar@{=}[d] &E_F  \ar[r]^p
\ar[d]^{\alpha^*} & Q \ar[r]\ar@{=}[d] & 1,&\;\;\;\;E_F
\ar[r]^\varepsilon&D \\
 0 \ar[r]& B \ar[r]^{j'}  &E_{F'}  \ar[r]^{p'}   & Q \ar[r]&
 1,&\;\;\;\;E_{F'}
\ar[r]^{\varepsilon'}&D}
\end{diagram}
\end{align*}
It remains to show that $\varepsilon'\alpha^*=\varepsilon$.

Since $\alpha:F\ri F'$ is a homotopy, and by (\ref{ht12}),
$Fu=x_{\psi(u)}=F'u$. Therefore
$x_{\psi(u)}=d(\alpha_u)x_{\psi(u)},$ or $d(\alpha_u)=1$. Then,
\[\begin{aligned}\varepsilon'\alpha^*(b,u)&=\varepsilon'(b+\alpha_u,u)=d(b+\alpha_u)x_{\psi(u)}\\
\;\;&
=d(b)d(\alpha_u)x_{\psi(u)}=d(b)x_{\psi(u)}=\varepsilon(b,u).\end{aligned}\]
That means $\mathcal E_F$ and $\mathcal E_{F'}$ are equivalent.

Conversely, if $\alpha^*:E\rightarrow E'$ is an isomorphism then
$$\alpha^*(b,u)=(b+\alpha_u,u),$$ where $\alpha:Q\ri B$ is a function
satisfying $\alpha_1=0$. By retracing our steps, $\alpha$ is a
homotopy of $F$ and $F'$.

{\it Step 2: $\Omega $ is a surjection.}

Assume that $\mathcal E$ is an extension $ E$ of $B$ by $Q$ of type
$B\rightarrow D$ inducing $\psi: Q\rightarrow \Coker d$ as in the
commutative diagram (\ref{bd10}).
We prove that $\mathcal E$ is equivalent to an extension $\mathcal
E_F$ which is associated to a Gr-functor
$(F,\widetilde{F}):\mathrm{Dis}Q\ri\mathbb P$.

 Let $\mathbb P'$ be the Gr-category associated to the crossed module
$(B,E,j,\theta')$. Then, by Proposition \ref{t1}, the pair
$(id_B,\varepsilon)$ determines a single Gr-functor
$(K,\widetilde{K}):\mathbb P'\ri\mathbb P$. Since $\pi_0\mathbb
P'=Q,\pi_1\mathbb P'=0$, then 
$S_{\mathbb P'}=\mathrm{Dis}Q.$

 Choose the stick $(e_u,i_e), e\in
E, u\in Q$ in $\mathbb P'$, i.e., $\{e_u\}$ is a representative of
$Q$ in $E$. By (\ref{ct1a}), the canonical Gr-functor
$(H',\widetilde{H}'): \Dis Q\ri \mathbb P'$ is given by
$$H'(u)=e_u, \widetilde{H}'_{u,v}=-i_{e_u+e_v}=h'(u,v).$$
Then the composition $F=KH'$ determines a Gr-functor $\Dis Q\ri
\mathbb P$, where
$$F(u)=\varepsilon(e_u),\widetilde{F}_{u,v}=K(\widetilde{H}'_{u,v})=h'(u,v).$$
By Lemma \ref{mrlk},  we determine an extension $\mathcal E_F$ of
the crossed product  $E_0=[B,\varphi,h',Q]$, which is associated to
$(F,\widetilde{F})$. We now prove that $\mathcal E$ and $\mathcal
E_F$ is equivalent, i.e., the following diagram commutes
\begin{align*}  \begin{diagram}
\xymatrix{\mathcal E_F:\;\;\; 0 \ar[r]& B \ar[r]^{j_0} \ar@{=}[d]
&E_0 \ar[r]^{p_0} \ar[d]^\eta & Q \ar[r]\ar@{=}[d] & 1,&\;\;\;\;E_0
\ar[r]^{\varepsilon_0}&D \\
\mathcal E:\;\;\; 0 \ar[r]& B \ar[r]^{j}  &E \ar[r]^{p}   & Q
\ar[r]& 1,&\;\;\;\;E \ar[r]^{\varepsilon}&D}
\end{diagram}
\end{align*}
The representative ${e_u}$ satisfies:
\begin{equation}\label{ct23}
\varphi(u)c=\mu_{e_u}(c), \ c\in B.
\end{equation}
\begin{equation}\label{ct24}
e_u+e_v=-i_{e_u+e_v}+e_{uv}=h'(u,v)+e_{uv}.
\end{equation}
The relation (\ref{ct23}) holds since $(id_B,\varepsilon)$ is a
morphism of crossed modules, the relation (\ref{ct24}) holds thanks
to the definition of morphism in $\mathbb P'$.

Each element of  $E$ is presented uniquely as $b+e_u, b\in B$, so
one can define a map
$$\eta:E_0\ri E, \ (b,u)\mapsto b+e_u$$
By the relations (\ref{ct23}) and (\ref{ct24}), $\eta$ is an
isomorphism.

 Finally, choose the representative $e_u$ such that
\begin{equation*}
\varepsilon(e_u)=x_{\psi(u)}.
\end{equation*}
Indeed, from  (\ref{bd10}) we have $q(\varepsilon(e_u))=\psi
p(e_u)=\psi(u).$ Then
$$\varepsilon \eta(b,u)=\varepsilon (b+e_u)=\varepsilon(b)\varepsilon(e_u)=d(b)x_{\psi(u)}=\varepsilon_0(b,u).$$
Thus, $\mathcal E$ and $\mathcal E_F$ are equivalent. This completes
the proof.\end{proof}

 Now, the second assertion of Theorem 5.2 \cite{Br94} is
obtained as follows. Firstly, there is a natural bijection
$$\Hom[\Dis Q, \mathbb P]\leftrightarrow \Hom[\Dis Q, S_{\mathbb P}]$$
Since $\pi_0(\Dis Q)=Q, \pi_1(S_{\mathbb P})=\Ker d$, then it
follows from Theorem \ref{schr} and Proposition \ref{md0} that
$$\mathrm{Ext}_{B\ri D}(Q, B,\psi)\leftrightarrow
H^2(Q,\mathrm{Ker }d).$$

\begin{center}
{}
\end{center}

\end{document}